\documentclass[english,11pt]{article}
\usepackage{amsfonts}
\usepackage[latin1]{inputenc}
\usepackage{amsmath,amsthm,amssymb}
\usepackage{graphicx}
\usepackage{color}

\def\neweq#1{\begin{equation}\label{#1}}
\def\endeq{\end{equation}}

\newtheorem{theorem}{Theorem}[section]

\input{tcilatex}
\begin{document}

\title{\textbf{Radial and nonradial solutions for a semilinear elliptic
system of Schrödinger type}\\
"to appear in FUNKCIALAJ EKVACIOJ"}
\author{Dragos-Patru Covei\break \\
{\small \ Constantin Brancusi University of \ Tg-Jiu and West University of
Timi\c{s}oara, Romania}\\
{\small E-mail: c\texttt{ovdra@yahoo.com}}}
\date{}
\maketitle

\begin{abstract}
In this article we consider the system of equations $\Delta
u_{i}=p_{i}\left( x\right) f_{i}\left( u_{1},...,u_{d}\right) $ for $%
i=1,...,d$ on $R^{N}$, $N\geq 3$ and $d\in \left\{ 1,2,3,4,...\right\} $. We
prove that the considered system has a bounded positive entire solution
under some conditions on $p_{i}$ and $f_{i}$. Also, we give a necessary
condition as well as a sufficient condition for a positive radial solution
to be large. The method of proving theorems is essentially based on a
successive approximation. Furthermore, a non-radially symmetric solution is
obtained by using a lower and upper solution method. \\[3pt]
\end{abstract}

\baselineskip16pt \renewcommand{\theequation}{\arabic{section}.%
\arabic{equation}} \catcode`@=11 \@addtoreset{equation}{section} \catcode%
`@=12

\textbf{2000 Mathematics Subject Classification}: 35J61, 35J91.

\textbf{Key words}: Entire solution; Large solution; Elliptic system.

\section{Introduction}

In this paper we study the existence of solutions for the semilinear
elliptic system 
\begin{equation}
\left\{ 
\begin{array}{l}
\Delta u_{1}=p_{1}\left( x\right) f_{1}\left( u_{1},...,u_{d}\right) \text{, 
} \\ 
... \\ 
\Delta u_{d}=p_{d}\left( x\right) f_{d}\left( u_{1},...,u_{d}\right) \text{,}%
\end{array}%
\right. \text{in }\mathbb{R}^{N}\text{,}  \label{11}
\end{equation}%
where $N\geq 3$, $d\geq 1$, 
and the functions $p_{j}$, $f_{j}$ ($j=1,\ldots ,d$) are supposed to satisfy
the following hypothesis:

\medskip \noindent (P1)\quad $p_{j}:\mathbb{R}^{N}\rightarrow \left[
0,\infty \right) $ are locally Hölder continuous functions of exponent $%
\alpha \in \left( 0,1\right) $;

\noindent (C1)\quad $f_{j}:\left[ 0,\infty \right) ^{d}\rightarrow \left[
0,\infty \right) $ are continuously differentiable in each variable, \ $%
f_{j}\left( 0,..,0\right) =0$, and $f_{j}\left( s_{1},...,s_{d}\right) >0$
if $s_{i}>0$ for some $i=1,\ldots ,d$;

\noindent (C2)\quad $f_{j}$ are increasing in each variable;

\noindent (C3)\quad $\int_{1}^{\infty }[F\left( s\right) ]^{-1/2}ds=\infty $
\ \ ($F\left( s\right) =\int_{0}^{s}\underset{i=1}{\overset{d}{\Sigma }}%
f_{i}\left( t,...,t\right) dt$).

Similar problems to those we are analyzing here are related to steady-state
reaction-diffusion, subsonic fluid flows, electric potentials of some bodies
and control theory.

As an example, the first motivation for studying the above problem stems
from the article \cite{LAS} where the reader observe that such problems
arise from the description of the basic stochastic control theory. The
controls are to be designed so that the state of the system is constrained
to some region. Finding optimal controls is then shown to be equivalent to
finding large solutions for a second order semilinear elliptic partial
differential equation. \ In terms of the dynamic programming approach, an
explosive solution of (\ref{11}) corresponds to a value function (or Bellman
function) associated to an infinite exit cost (see \cite{LAS}).

Another motivation comes from the work of \cite{B} where the parabolic
problem corresponding to system (\ref{11}) \ are models of steady state of
non-linear heat conduction through a 2-components mixture.

The numerous applications that lead to favorable answers in order to
establish new and significant results for problem (\ref{11}) as well as the
recent results in the field motivate the study of more generally class of
problems (\ref{11}).

The main results of this paper are the following two theorems. They
substantially solve the open problem proposed in \cite{LAIR2} and complete
the results of \cite{CD2} where only sufficient conditions are obtained. 

\begin{theorem}
\label{1} Suppose that \textrm{(P1)}, \textrm{(C1)}--\textrm{(C3)} are
satisfied. If there exists a positive number $\varepsilon $ such that 
\begin{equation}
\int_{0}^{\infty }t^{1+\varepsilon }\underset{j=1}{\overset{d}{\Sigma }}%
\varphi _{j}\left( t\right) dt<\infty\quad\text{ where }\quad\varphi _{j}
\left( t\right) =\max_{\left\vert x\right\vert =t}p_{j}(x),  \label{5}
\end{equation}%
and $r^{2N-2}\underset{j=1}{\overset{d}{\Sigma }}\varphi _{j}\left( r\right) 
$ is nondecreasing for large $r$, then system \textrm{(\ref{11})} has a
nonnegative nontrivial bounded solution on $\mathbb{R}^{N}$. If, on the
other hand, $p_{j}$ satisfy%
\begin{equation}
\int_{0}^{\infty }t\underset{j=1}{\overset{d}{\Sigma }}\psi _{j}\left(
t\right) dt=\infty \quad\text{ where }\quad\psi _{j}\left( t\right)
=\min_{\left\vert x\right\vert =t}p_{j}\left( x\right) ,  \label{5b}
\end{equation}%
and $r^{2N-2}\underset{j=1}{\overset{d}{\Sigma }}\psi _{j}\left( r\right) $
is nondecreasing for large $r$, then system \textrm{(\ref{11})} has no
nonnegative nontrivial entire bounded radial solution on $\mathbb{R}^{N}$.
\end{theorem}

\begin{theorem}
\label{2}\textit{Assume that }$p_{j}$\textit{\ }$:\mathbb{R}^{N}\rightarrow %
\left[ 0,\infty \right) $ $(j=1,\ldots ,d)$ \textit{\ are spherically
symmetric continuous functions $($i.e. }$p_{j}\left( x\right) =p_{j}\left(
\left\vert x\right\vert \right) $\textit{$)$. If }\ $f_{j}$ $(j=1,\ldots ,d)$
\textit{\ satisfy \textrm{(C1)} --\textrm{(C3)}, then the problem \textrm{(%
\ref{11})} has a nonnegative nontrivial entire radial solution. Suppose
furthermore that }$r^{2N-2}\underset{j=1}{\overset{d}{\Sigma }}p_{j}\left(
r\right) $\textit{\ is nondecreasing for large }$r$\textit{.} If $p_{j}$
satisfies 
\begin{equation}
\int_{0}^{\infty }\frac{1}{t^{N-1}}\int_{0}^{t}s^{N-1}p_{j}\left( s\right)
dsdt=\infty \quad \text{ for all }\quad j=1,\ldots ,d,\text{ }  \label{12}
\end{equation}%
then any nonnegative nontrivial solution $\left( u_{1},...,u_{d}\right) $ of 
\textrm{(\ref{11})} is large. Conversely, if \textrm{(\ref{11})} has a
nonnegative entire large solution, then $p_{j}$ satisfy%
\begin{equation}
\int_{0}^{\infty }r^{1+\varepsilon }\underset{j=1}{\overset{d}{\Sigma }}%
p_{j}\left( r\right) dr=\infty ,  \label{13}
\end{equation}%
for every $\varepsilon >0$.
\end{theorem}

\section{Preliminary result}

We recall the following definition of lower and upper solution which are our
main tool in the proof of the solvability of problem (\ref{11}).

\begin{definition}
A function $\left( w_{1},...,w_{d}\right) \in \left[ C_{loc}^{2,\alpha
}\left( \mathbb{R}^{N}\right) \right] ^{d}$ ($\alpha \in \left( 0,1\right) $%
) is called a lower solution of the problem (\ref{11}) if%
\begin{equation*}
\Delta w_{i}\geq p_{i}\left( x\right) \ f_{i}\left( w_{1},...,w_{d}\right) 
\text{ in }\mathbb{R}^{N}\text{ for all \ }i=1,...,d\text{.}
\end{equation*}%
$\ $
\end{definition}

\begin{definition}
We say that $\left( v_{1},...,v_{d}\right) \in \left[ C_{loc}^{2,\alpha
}\left( \mathbb{R}^{N}\right) \right] ^{d}(\alpha \in \left( 0,1\right) )$
is an upper solution of the problem (\ref{11}) if%
\begin{equation*}
\Delta v_{i}\leq p_{i}\left( x\right) \ f_{i}\left( v_{1},...,v_{d}\right) 
\text{ in }\mathbb{R}^{N}\text{ for all\ }i=1,...,d\text{.}
\end{equation*}
\end{definition}

We need the following lemma which can be found in \cite[Theorem 5.1, pp. 146]%
{NK}:

\begin{lemma}
\label{lu} Make the same assumptions on $p_{j}$ and \ $f_{j}$ ($j=1,...,d$)
as in Theorem \ref{1}. If the problem (\ref{11}) has a pair of upper and
lower bounded solutions $\left( v_{1},...,v_{d}\right) $ and $\left(
w_{1},...,w_{d}\right) $ fulfilling $w_{i}(x)\leq v_{i}(x),$ $%
i=1,...,d,\forall x\in \mathbb{R}^{N}$ then there exists a bounded function $%
\left( u_{1},...,u_{d}\right) $ belonging to $\left[ C_{loc}^{2,\alpha
}\left( \mathbb{R}^{N}\right) \right] ^{d}$ ($\alpha \in \left( 0,1\right) $%
) with%
\begin{equation*}
w_{i}(x)\leq u_{i}(x)\leq v_{i}(x)\text{, }i=1,...,d\text{, }\forall x\in 
\mathbb{R}^{N}
\end{equation*}%
and satisfying (\ref{11}).
\end{lemma}

\section{Proofs of main results}

In this section we give the proofs of Theorems \ref{1} and \ref{2}.

\subparagraph{Proof of the Theorem \protect\ref{1}}

Assume that (\ref{5}) holds. In this section, we use the method of upper and
lower solutions to show the existence of positive bounded solutions of the
problem (\ref{11}). Our aim is to construct an upper solution $\left(
v_{1},...,v_{d}\right) $ and a lower solution $\left( w_{1},...,w_{d}\right) 
$ \ for the problem (\ref{11}). On this purpose, we first prove the
existence of $\left( w_{1},...,w_{d}\right) $ to the system%
\begin{equation}
\left\{ 
\begin{array}{l}
\Delta w_{1}\left( r\right) =\varphi _{1}\left( r\right) f_{1}\left(
w_{1},...,w_{d}\right) \text{ for }r:=\left\vert x\right\vert \text{, } \\ 
... \\ 
\Delta w_{d}\left( r\right) =\varphi _{d}\left( r\right) f_{d}\left(
w_{1},...,w_{d}\right) \text{ for }r:=\left\vert x\right\vert \text{, }%
\end{array}%
\right. \text{in }\mathbb{R}^{N}.  \label{6}
\end{equation}%
Observe that we can rewrite (\ref{6}) as follows: 
\begin{equation*}
\left\{ 
\begin{array}{l}
\left( r^{N-1}w_{1}\left( r\right) \right) ^{\prime }=r^{N-1}\varphi
_{1}\left( r\right) f_{1}\left( w_{1},...,w_{d}\right) , \\ 
... \\ 
\left( r^{N-1}w_{d}\left( r\right) \right) ^{\prime }=r^{N-1}\varphi
_{d}\left( r\right) f_{d}\left( w_{1},...,w_{d}\right) .%
\end{array}%
\right.
\end{equation*}%
Then radial solutions of (\ref{6}) are any solution $\left(
w_{1},...,w_{d}\right) $ of the integral equations%
\begin{equation*}
\left\{ 
\begin{array}{l}
w_{1}\left( r\right) =\frac{1}{d}+\int_{0}^{r}\frac{1}{t^{N-1}}%
\int_{0}^{t}s^{N-1}\varphi _{1}\left( s\right) f_{1}\left( w_{1}\left(
s\right) ,...,w_{d}\left( s\right) \right) dsdt,\text{ } \\ 
... \\ 
w_{d}\left( r\right) =\frac{1}{d}+\int_{0}^{r}\frac{1}{t^{N-1}}%
\int_{0}^{t}s^{N-1}\varphi _{d}\left( s\right) f_{d}\left( w_{1}\left(
s\right) ,...,w_{d}\left( s\right) \right) dsdt.%
\end{array}%
\right.
\end{equation*}%
To establish a solution to this system, we use successive approximation.
Define sequences $\left\{ w_{j}^{k}\right\} _{j=1,...,d}^{k\geq 1}$ on $%
\left[ 0,\infty \right) $ by%
\begin{equation*}
\left\{ 
\begin{array}{l}
w_{1}^{0}=...=w_{d}^{0}=\frac{1}{d},\text{ }r\geq 0, \\ 
w_{i}^{k}\left( r\right) =\frac{1}{d}+\int_{0}^{r}\frac{1}{t^{N-1}}%
\int_{0}^{t}s^{N-1}\varphi _{i}\left( s\right) f_{i}\left( w_{1}^{k-1}\left(
s\right) ,...,w_{d}^{k-1}\left( s\right) \right) dsdt,\text{ }i=1,...,d.%
\end{array}%
\right.
\end{equation*}%
We remark that, for all $r\geq 0,$ $j=1,...,d$ and $k\in N$ 
\begin{equation*}
w_{j}^{k}\left( r\right) \geq \frac{1}{d}\text{.}
\end{equation*}%
Moreover, proceeding by induction we conclude $\left\{ w_{j}^{k}\right\}
_{j=1,...,d}^{k\geq 1}$ are non-decreasing sequence on $\left[ 0,\infty
\right) $. We note that $\left\{ w_{j}^{k}\right\} _{j=1,...,d}^{k\geq 1}$
satisfy%
\begin{equation*}
\left\{ 
\begin{array}{l}
\left[ r^{N-1}\left( w_{1}^{k}\right) ^{\prime }\right] ^{\prime
}=r^{N-1}\varphi _{1}\left( r\right) f_{1}\left(
w_{1}^{k-1},...,w_{d}^{k-1}\right) ,\text{ } \\ 
... \\ 
\left[ r^{N-1}\left( w_{d}^{k}\right) ^{\prime }\right] ^{\prime
}=r^{N-1}\varphi _{d}\left( r\right) f_{d}\left(
w_{1}^{k-1},...,w_{d}^{k-1}\right) .%
\end{array}%
\right.
\end{equation*}%
By the monotonicity of $\left\{ w_{j}^{k}\right\} _{j=1,...,d}^{k\geq 1}$ we
have the inequalities%
\begin{equation}
\left\{ 
\begin{array}{l}
\left[ r^{N-1}\left( w_{1}^{k}\right) ^{\prime }\right] ^{\prime
}=r^{N-1}\varphi _{1}\left( r\right) f_{1}\left(
w_{1}^{k-1},...,w_{d}^{k-1}\right) \leq r^{N-1}\varphi _{1}\left( r\right) 
\overset{d}{\underset{i=1}{\Sigma }}f_{i}\left( \overset{d}{\underset{j=1}{%
\Sigma }}w_{j}^{k}\left( r\right) ,...,\overset{d}{\underset{j=1}{\Sigma }}%
w_{j}^{k}\left( r\right) \right) , \\ 
... \\ 
\left[ r^{N-1}\left( w_{d}^{k}\right) ^{\prime }\right] ^{\prime
}=r^{N-1}\varphi _{d}\left( r\right) f_{d}\left(
w_{1}^{k-1},...,w_{d}^{k-1}\right) \leq r^{N-1}\varphi _{d}\left( r\right) 
\overset{d}{\underset{i=1}{\Sigma }}f_{i}\left( \overset{d}{\underset{j=1}{%
\Sigma }}w_{j}^{k}\left( r\right) ,...,\overset{d}{\underset{j=1}{\Sigma }}%
w_{j}^{k}\left( r\right) \right) .%
\end{array}%
\right.  \label{8}
\end{equation}%
Thus, summing up inequalities, we have%
\begin{equation}
\left[ r^{N-1}\overset{d}{\underset{i=1}{\Sigma }}\left( w_{i}^{k}\right)
^{\prime }\right] ^{\prime }\leq \text{ }r^{N-1}\overset{d}{\underset{i=1}{%
\Sigma }}\varphi _{i}\left( r\right) \overset{d}{\underset{i=1}{\Sigma }}%
f_{i}\left( \overset{d}{\underset{j=1}{\Sigma }}w_{j}^{k}\left( r\right)
,...,\overset{d}{\underset{j=1}{\Sigma }}w_{j}^{k}\left( r\right) \right) .
\label{88}
\end{equation}%
Choose $R>0$ so that $r^{2N-2}\underset{j=1}{\overset{d}{\Sigma }}\varphi
_{j}\left( r\right) $ are non-decreasing for $r\geq R$. We are now ready to
show that $w_{j}^{k}\left( R\right) $ and $\left( w_{j}^{k}\left( R\right)
\right) ^{\prime }$, both of which are nonnegative, are bounded above
independent of $k$. To do this, let 
\begin{equation*}
\phi _{j}^{R}=\max \{\varphi _{j}\left( r\right) :0\leq r\leq R\}\text{, }%
j=1,...,d\text{.}
\end{equation*}%
Using this and the fact that $\left( w_{j}^{k}\right) ^{\prime }\geq 0$, we
note that (\ref{8}) yields%
\begin{equation*}
\left\{ 
\begin{array}{l}
\left[ \left( w_{1}^{k}\right) ^{\prime }\right] ^{\prime }\leq \phi _{1}^{R}%
\overset{d}{\underset{i=1}{\Sigma }}f_{i}\left( \overset{d}{\underset{i=1}{%
\Sigma }}w_{i}^{k},...,\overset{d}{\underset{i=1}{\Sigma }}w_{i}^{k}\right) ,
\\ 
... \\ 
\left[ \left( w_{d}^{k}\right) ^{\prime }\right] ^{\prime }\leq \phi _{d}^{R}%
\overset{d}{\underset{i=1}{\Sigma }}f_{i}\left( \overset{d}{\underset{i=1}{%
\Sigma }}w_{i}^{k},...,\overset{d}{\underset{i=1}{\Sigma }}w_{i}^{k}\right) 
\text{. }%
\end{array}%
\right.
\end{equation*}%
This implies that 
\begin{equation*}
\underset{i=1}{\overset{d}{\Sigma }}\left[ \left( w_{i}^{k}\right) ^{\prime }%
\right] ^{\prime }\leq \left( \underset{i=1}{\overset{d}{\Sigma }}\phi
_{i}^{R}\right) \overset{d}{\underset{i=1}{\Sigma }}f_{i}\left( \overset{d}{%
\underset{i=1}{\Sigma }}w_{i}^{k},...,\overset{d}{\underset{i=1}{\Sigma }}%
w_{i}^{k}\right) .
\end{equation*}%
On the other hand, multiplying this equation by $\left( \underset{i=1}{%
\overset{d}{\Sigma }}w_{i}^{k}\right) ^{\prime }$ and integrating we also
have%
\begin{equation}
\left[ \left( \underset{i=1}{\overset{d}{\Sigma }}w_{i}^{k}\left( r\right)
\right) ^{\prime }\right] ^{2}\leq 2\left( \underset{i=1}{\overset{d}{\Sigma 
}}\phi _{i}^{R}\right) \int_{1}^{\overset{d}{\underset{i=1}{\Sigma }}%
w_{i}^{k}\left( r\right) }\overset{d}{\underset{i=1}{\Sigma }}f_{i}\left(
s,...,s\right) ds\text{, }0\leq r\leq R.  \label{9}
\end{equation}%
Integrating the above inequality, we see that 
\begin{equation*}
\int_{1}^{\overset{d}{\underset{i=1}{\Sigma }}w_{i}^{k}\left( R\right) }%
\left[ \int_{1}^{t}\overset{d}{\underset{i=1}{\Sigma }}f_{i}\left(
s,...,s\right) ds\right] ^{-1/2}dt\leq \sqrt{2\left( \underset{i=1}{\overset{%
d}{\Sigma }}\phi _{i}^{R}\right) }R\text{.}
\end{equation*}%
It follows from the above relation and by the assumption C3) that $\overset{d%
}{\underset{i=1}{\Sigma }}w_{i}^{k}\left( R\right) $ is bounded above
independent of $k$. Using this fact in (\ref{9}) shows that the same is true
of $\left( \overset{d}{\underset{i=1}{\Sigma }}w_{i}^{k}\left( R\right)
\right) ^{\prime }$. Thus, the sequences $w_{i}^{k}\left( R\right) $ and $%
\left( w_{i}^{k}\left( R\right) \right) ^{\prime }$ are bounded above
independent of $k$. Finally, we show that the non-decreasing sequences $%
w_{i}^{k}$ is bounded for all $r\geq 0$ and all $k$. Multiplying (\ref{88})
by $r^{N-1}\left( \overset{d}{\underset{i=1}{\Sigma }}w_{i}^{k}\right)
^{\prime }$ and integrating gives%
\begin{eqnarray*}
&&\left[ r^{N-1}\left( \overset{d}{\underset{i=1}{\Sigma }}w_{i}^{k}\left(
r\right) \right) ^{\prime }\right] ^{2} \\
&\leq &\left[ R^{N-1}\left( \overset{d}{\underset{i=1}{\Sigma }}%
w_{i}^{k}\left( R\right) \right) ^{\prime }\right] ^{2}+2\int_{R}^{r}s^{2N-2}%
\overset{d}{\underset{i=1}{\Sigma }}\varphi _{i}\left( s\right) \overset{d}{%
\underset{i=1}{\Sigma }}f_{i}\left( \overset{d}{\underset{i=1}{\Sigma }}%
w_{i}^{k},...,\overset{d}{\underset{i=1}{\Sigma }}w_{i}^{k}\right) \left( 
\overset{d}{\underset{i=1}{\Sigma }}w_{i}^{k}\right) ^{\prime }ds\text{,}
\end{eqnarray*}%
for $r\geq R$. Noting that, by the monotonicity of $s^{2N-2}\overset{d}{%
\underset{i=1}{\Sigma }}\varphi _{i}\left( s\right) $ for $s\geq R$, we get 
\begin{equation*}
\left[ r^{N-1}\left( \overset{d}{\underset{i=1}{\Sigma }}w_{i}^{k}\left(
r\right) \right) ^{\prime }\right] ^{2}\leq C+2r^{2N-2}\overset{d}{\underset{%
i=1}{\Sigma }}\varphi _{i}\left( r\right) F\left( \overset{d}{\underset{i=1}{%
\Sigma }}w_{i}^{k}\left( r\right) \right) ,
\end{equation*}%
where $C=\left[ R^{N-1}\left( \overset{d}{\underset{i=1}{\Sigma }}%
w_{i}^{k}\left( R\right) \right) ^{\prime }\right] ^{2}$, which yields%
\begin{equation}
\left( \overset{d}{\underset{i=1}{\Sigma }}w_{i}^{k}\right) ^{\prime }\leq 
\sqrt{C}r^{1-N}+\sqrt{2\overset{d}{\underset{i=1}{\Sigma }}\varphi
_{i}\left( r\right) }\left[ F\left( \overset{d}{\underset{i=1}{\Sigma }}%
w_{i}^{k}\left( r\right) \right) \right] ^{1/2},  \label{10}
\end{equation}%
and hence%
\begin{equation*}
\frac{d}{dr}\int_{\overset{d}{\underset{i=1}{\Sigma }}w_{i}^{k}\left(
R\right) }^{\overset{d}{\underset{i=1}{\Sigma }}w_{i}^{k}\left( r\right) }%
\left[ F\left( t\right) \right] ^{-1/2}dt\leq \sqrt[2]{C}r^{1-N}\left[
F\left( \overset{d}{\underset{i=1}{\Sigma }}w_{i}^{k}\left( r\right) \right) %
\right] ^{-1/2}+\left( 2\overset{d}{\underset{i=1}{\Sigma }}\varphi
_{i}\left( r\right) \right) ^{1/2}.
\end{equation*}%
Integrating this and using the fact that%
\begin{eqnarray*}
\left( 2\overset{d}{\underset{i=1}{\Sigma }}\varphi _{i}\left( s\right)
\right) ^{1/2} &=&\left( 2s^{1+\varepsilon }\overset{d}{\underset{i=1}{%
\Sigma }}\varphi _{i}\left( s\right) s^{-1-\varepsilon }\right) ^{1/2} \\
&\leq &s^{1+\varepsilon }\left( \overset{d}{\underset{i=1}{\Sigma }}\varphi
_{i}\left( s\right) \right) +s^{-1-\varepsilon }
\end{eqnarray*}%
for each $\varepsilon >0$, we have%
\begin{equation}
\begin{array}{ll}
\int_{\overset{d}{\underset{i=1}{\Sigma }}w_{i}^{k}\left( R\right) }^{%
\overset{d}{\underset{i=1}{\Sigma }}w_{i}^{k}\left( r\right) }\left[ F\left(
t\right) \right] ^{-1/2}dt & \leq \sqrt[2]{C}\int_{R}^{r}t^{1-N}\left[
F\left( \overset{d}{\underset{i=1}{\Sigma }}w_{i}^{k}\left( t\right) \right) %
\right] ^{-1/2}dt \\ 
& +\left[ \int_{R}^{r}t^{1+\varepsilon }\overset{d}{\underset{i=1}{\Sigma }}%
\varphi _{i}\left( t\right) dt+\int_{R}^{r}t^{-1-\varepsilon }dt\right] \\ 
& \leq \sqrt[2]{C}\left[ F\left( \overset{d}{\underset{i=1}{\Sigma }}%
w_{i}^{k}\left( R\right) \right) \right] ^{-1/2}\int_{R}^{r}t^{1-N}dt \\ 
& +\left[ \int_{R}^{r}t^{1+\varepsilon }\left( \overset{d}{\underset{i=1}{%
\Sigma }}\varphi _{i}\left( t\right) \right) dt+\frac{1}{\varepsilon
R^{\varepsilon }}\right] \text{.}%
\end{array}
\label{111}
\end{equation}%
The above relation is needed in proving the bounded of the function $\left\{
w_{j}^{k}\right\} _{j=1,...,d}^{k\geq 1}$ in the following. Indeed, since
for each $\varepsilon >0$ the right side of this inequality is bounded
independent of $k$ (note that $w_{i}^{k}\left( t\right) \geq 1/d$ for all $%
i=1,...,d$), so is the left side and hence, in light of C3), the sequence $%
\left\{ \overset{d}{\underset{i=1}{\Sigma }}w_{i}^{k}\right\} ^{k\geq 1}$ is
a bounded sequence and so $\left\{ w_{j}^{k}\right\} _{j=1,...,d}^{k\geq 1}$
are bounded sequence. Thus $\left\{ w_{j}^{k}\right\} _{j=1,...,d}^{k\geq
1}\rightarrow \left\{ w_{j}\right\} _{j=1,...,d}$ as $k\rightarrow \infty $
and the limit functions $\left\{ w_{j}\right\} _{j=1,...,d}$ are positive
entire solutions of system (\ref{6}). It is obvious that every solution of (%
\ref{6}) is a lower solution to (\ref{11}).

We let $M$ be the least upper bound of $\overset{d}{\underset{i=1}{\Sigma }}%
w_{i}$ and note that 
\begin{equation*}
M=\lim_{r\rightarrow \infty }\overset{d}{\underset{i=1}{\Sigma }}w_{i}\left(
r\right) \text{.}
\end{equation*}%
Now let $\psi _{i}\left( t\right) =\min_{\left\vert x\right\vert
=t}p_{i}\left( x\right) $ ($i=1,...,d$)\ and $\left\{ v_{j}\right\}
_{j=1,...,d}$ be the positive increasing bounded solutions of 
\begin{equation*}
\left\{ 
\begin{array}{l}
v_{1}\left( r\right) =M+\int_{0}^{r}\frac{1}{t^{N-1}}\int_{0}^{t}s^{N-1}\psi
_{1}\left( s\right) f_{1}\left( v_{1}\left( s\right) ,...,v_{d}\left(
s\right) \right) dsdt,\text{ } \\ 
... \\ 
v_{d}\left( r\right) =M+\int_{0}^{r}\frac{1}{t^{N-1}}\int_{0}^{t}s^{N-1}\psi
_{d}\left( s\right) f_{d}\left( v_{1}\left( s\right) ,...,v_{d}\left(
s\right) \right) dsdt,%
\end{array}%
\right.
\end{equation*}%
which, of course, satisfies (\ref{6}) with $w_{j}$ ($j=1,...,d$) replaced
with $v_{j}$ ($j=1,...,d$) and $\varphi _{j}$ ($j=1,...,d$) replaced with $%
\psi _{j}$ ($j=1,...,d$). It is also clear that $v_{j}\geq M$. If we modify
the proof of the lower solution $\left( w_{1},...,w_{d}\right) $ we obtain
the existence of an upper solution $\left( v_{1},...,v_{d}\right) $. Hence
the standard upper-lower solution principle implies that (\ref{11}) has a
solution $\left( u_{1},...,u_{d}\right) $ such that%
\begin{equation*}
w_{i}(x)\leq u_{i}(x)\leq v_{i}(x),\text{ }i=1,...,d,
\end{equation*}%
(see Lemma \ref{lu}).

Now assume that (\ref{5b}) holds. Assume to the contrary that there exist a
nonnegative nontrivial entire bounded radial solution $\left(
u_{1},...,u_{d}\right) $ on $\mathbb{R}^{N}$ for the system (\ref{11}).
Assuming $M_{i}=\sup_{x\in \mathbb{R}^{N}}u_{i}\left( x\right) $ ($i=1,...,d$%
)\ and knowing that $u_{i}^{\prime }\geq 0$, we get $\lim_{r\rightarrow
\infty }u_{i}=M_{i}$. Thus there exists $R>0$ such that $u_{i}\geq M_{i}/2$ (%
$i=1,...,d$)\ for $r\geq R$. From conditions of $f_{i}$, it follows that%
\begin{equation*}
f_{i}\left( u_{1},...,u_{d}\right) \geq f_{i}\left(
M_{1}/2,...,M_{d}/2\right) :=c_{0}^{i}\text{, for }r\geq R\text{ and }%
i=1,...,d\text{.}
\end{equation*}%
Using this we get the following%
\begin{eqnarray}
u_{1} &\geq &u_{1}\left( 0\right)
+\int_{0}^{r}t^{1-N}\int_{0}^{t}s^{N-1}\psi _{1}\left( s\right) f_{1}\left(
u_{1}\left( s\right) ,...,u_{d}\left( s\right) \right) dsdt  \notag \\
&&...  \label{sis4} \\
u_{d} &\geq &u_{d}\left( 0\right)
+\int_{0}^{r}t^{1-N}\int_{0}^{t}s^{N-1}\psi _{d}\left( s\right) f_{d}\left(
u_{1}\left( s\right) ,...,u_{d}\left( s\right) \right) dsdt.  \notag
\end{eqnarray}%
Rearranging the terms, and by using these conditions in (\ref{sis4}) follows%
\begin{equation*}
\overset{d}{\underset{i=1}{\Sigma }}u_{i}\left( r\right) \geq dm_{1}+m_{2}%
\frac{1}{N}\int_{R}^{r}t\overset{d}{\underset{i=1}{\Sigma }}\psi _{i}\left(
t\right) dt\rightarrow \infty \text{ as }r\rightarrow \infty ,\text{ }
\end{equation*}%
where $m_{1}=\min \{u_{1}\left( 0\right) ,...,u_{d}\left( 0\right) \}$ and $%
m_{2}:=\min \{c_{0}^{1},...,c_{0}^{d}\}$. A contradiction to the boundedness
of $\overset{d}{\underset{i=1}{\Sigma }}u_{i}\left( r\right) $. This proves
Theorem \ref{1}.

\subparagraph{Proof of the Theorem \protect\ref{2}.}

It is known from \cite[Theorem 2]{LAIR,Y} that the problem%
\begin{equation}
\Delta z\left( r\right) =\overset{d}{\underset{i=1}{\Sigma }}p_{i}\left(
r\right) \overset{d}{\underset{i=1}{\Sigma }}f_{i}\left( z\left( r\right)
,...,z\left( r\right) \right) \text{ for }r:=\left\vert x\right\vert \text{, 
}x\in \mathbb{R}^{N}  \label{prob}
\end{equation}%
has a non-negative non-trivial entire solution. Moreover, for each $R>0$,
there exists $c_{R}>0$ such that $z\left( R\right) \leq c_{R}$. Due to the
fact that $z$ is radial, we have%
\begin{equation*}
z\left( r\right) =z\left( 0\right) +\int_{0}^{r}\frac{1}{t^{N-1}}%
\int_{0}^{t}s^{N-1}\overset{d}{\underset{i=1}{\Sigma }}p_{i}\left( s\right) 
\overset{d}{\underset{i=1}{\Sigma }}f_{i}\left( z\left( s\right)
,...,z\left( s\right) \right) dsdt\text{ for all }r\geq 0.
\end{equation*}%
We choose $\beta _{1}\in \left( 0,z\left( 0\right) \right] $. Define the
sequences $\left\{ u_{j}^{k}\right\} _{j=1,...,d}^{k\geq 1}$ on $\left[
0,\infty \right) $ by%
\begin{equation*}
\left\{ 
\begin{array}{l}
u_{1}^{0}=...=u_{d}^{0}=\beta _{1}\text{ for all }r\geq 0 \\ 
u_{1}^{k}\left( r\right) =\beta _{1}+\int_{0}^{r}\frac{1}{t^{N-1}}%
\int_{0}^{t}s^{N-1}p_{1}\left( s\right) f_{1}\left( u_{1}^{k-1}\left(
s\right) ,...,u_{d}^{k-1}\left( s\right) \right) dsdt, \\ 
... \\ 
u_{d}^{k}\left( r\right) =\beta _{1}+\int_{0}^{r}\frac{1}{t^{N-1}}%
\int_{0}^{t}s^{N-1}p_{d}\left( s\right) f_{d}\left( u_{1}^{k-1}\left(
s\right) ,...,u_{d}^{k-1}\left( s\right) \right) dsdt,\text{ }d\geq 1.%
\end{array}%
\right.
\end{equation*}%
With the same arguments as in the proof of Theorem \ref{1} we obtain
that\bigskip\ $\left\{ u_{j}^{k}\right\} _{j=1,...,d}^{k\geq 1}$ are
non-decreasing sequence on $\left[ 0,\infty \right) $. Because $z^{\prime
}\left( r\right) \geq 0$ follows $0<\beta _{1}\leq z\left( 0\right) \leq
z\left( r\right) $ for all $r\geq 0$ and so%
\begin{eqnarray*}
u_{1}^{1}\left( r\right) &=&\beta _{1}+\int_{0}^{r}\frac{1}{t^{N-1}}%
\int_{0}^{t}s^{N-1}p_{1}\left( s\right) f_{i}\left( u_{1}^{0}\left( s\right)
,...,u_{d}^{0}\left( s\right) \right) dsdt \\
&\leq &z\left( 0\right) +\int_{0}^{r}\frac{1}{t^{N-1}}\int_{0}^{t}s^{N-1}%
\overset{d}{\underset{i=1}{\Sigma }}p_{i}\left( s\right) \overset{d}{%
\underset{i=1}{\Sigma }}f_{i}\left( z\left( s\right) ,...,z\left( s\right)
\right) dsdt=z\left( r\right) .
\end{eqnarray*}%
Thus $u_{1}^{1}\left( r\right) \leq z\left( r\right) $. Similar arguments
show that%
\begin{equation*}
u_{j}^{k}\leq z\left( r\right) \text{ (}j=1,...,d\text{)\ for all }r\in %
\left[ 0,\infty \right) \text{ and }k\geq 1\text{.}
\end{equation*}%
Thus, for every $x\in \mathbb{R}^{N}$ the sequence $u_{j}^{k}\left(
\left\vert x\right\vert \right) $ ($j=1,...,d$)\ have subsequences, denoted
again by $u_{j}^{k}\left( \left\vert x\right\vert \right) $ ($j=1,...,d$),
converging and 
\begin{equation*}
\left( u_{1}\left( \left\vert x\right\vert \right) ,...,u_{d}\left(
\left\vert x\right\vert \right) \right) :=\left( \underset{k\rightarrow
\infty }{\lim }u_{1}^{k}\left( \left\vert x\right\vert \right) ,...,\underset%
{k\rightarrow \infty }{\lim }u_{d}^{k}\left( \left\vert x\right\vert \right)
\right) ,
\end{equation*}%
is an entire radial solution of system (\ref{11}).

Let $\left( u_{1},...,u_{d}\right) $ be any non-negative non-trivial entire
radial solution of (\ref{11}) and suppose that $p_{j}$ ($j=1,...,d$)
satisfies (\ref{12}). Since $u_{i}$ ($i=1,...,d$) is nontrivial and
non-negative, there exists $R>0$ so that $u_{i}\left( R\right) >0$. Since $%
u_{i}^{\prime }\geq 0$, we get $u_{i}\left( r\right) \geq u_{i}\left(
R\right) $ for $r\geq R$ and thus from%
\begin{equation*}
\left\{ 
\begin{array}{l}
u_{1}\left( r\right) =u_{1}\left( 0\right) +\int_{0}^{r}\frac{1}{t^{N-1}}%
\int_{0}^{t}s^{N-1}p_{1}\left( s\right) f_{1}\left( u_{1}\left( s\right)
,...,u_{d}\left( s\right) \right) dsdt,\text{ } \\ 
... \\ 
u_{d}\left( r\right) =u_{d}\left( 0\right) +\int_{0}^{r}\frac{1}{t^{N-1}}%
\int_{0}^{t}s^{N-1}p_{d}\left( s\right) f_{d}\left( u_{1}\left( s\right)
,...,u_{d}\left( s\right) \right) dsdt\text{, }%
\end{array}%
\right.
\end{equation*}%
we get 
\begin{equation*}
\left\{ 
\begin{array}{l}
u_{i}\left( r\right) =u_{i}\left( 0\right) +\int_{0}^{r}\frac{1}{t^{N-1}}%
\int_{0}^{t}s^{N-1}p_{i}\left( s\right) f_{i}\left( u_{1}\left( s\right)
,...,u_{d}\left( s\right) \right) dsdt\text{, \ }i=1,...,d\text{ } \\ 
\text{ \ \ \ \ \ \ }\geq u_{i}\left( R\right) +f_{i}\left( u_{1}\left(
R\right) ,...,u_{d}\left( R\right) \right) \int_{R}^{r}\frac{1}{t^{N-1}}%
\int_{R}^{t}s^{N-1}p_{i}\left( s\right) dsdt\rightarrow \infty \text{ as }%
r\rightarrow \infty .%
\end{array}%
\right.
\end{equation*}%
Conversely, if \ $f_{i}$ ($i=1,...,d$) satisfy (C1)-(C3) and $\left(
w_{1},..,w_{d}\right) $ is a nonnegative entire large solution of (\ref{11}%
), then $w_{i}$ satisfy%
\begin{equation*}
\left\{ 
\begin{array}{l}
\left[ r^{N-1}\left( w_{1}\right) ^{\prime }\right] ^{\prime }=p_{1}\left(
x\right) f_{1}\left( w_{1},...,w_{d}\right) , \\ 
... \\ 
\left[ r^{N-1}\left( w_{d}\right) ^{\prime }\right] ^{\prime }=p_{d}\left(
x\right) f_{d}\left( w_{1},...,w_{d}\right) .%
\end{array}%
\right.
\end{equation*}%
Then, using the monotonicity of \textit{\ }$r^{2N-2}\underset{j=1}{\overset{d%
}{\Sigma }}p_{j}\left( r\right) $\textit{\ } we can apply similar arguments
used in obtaining Theorem \ref{1} to get%
\begin{equation*}
\left( \overset{d}{\underset{i=1}{\Sigma }}w_{i}^{k}\left( r\right) \right)
^{\prime }\leq \sqrt{C}r^{1-N}+\sqrt{2\overset{d}{\underset{i=1}{\Sigma }}%
\varphi _{i}\left( r\right) }\sqrt{F\left( \overset{d}{\underset{i=1}{\Sigma 
}}w_{i}^{k}\left( r\right) \right) }.
\end{equation*}%
This conclusion and relation (\ref{111}) imply that%
\begin{equation}
\begin{array}{ll}
\int_{\overset{d}{\underset{i=1}{\Sigma }}w_{i}\left( R\right) }^{\overset{d}%
{\underset{i=1}{\Sigma }}w_{i}\left( r\right) }\left[ F\left( t\right) %
\right] ^{-1/2}dt & \leq \frac{\sqrt[2]{C}}{\left[ F\left( \overset{d}{%
\underset{i=1}{\Sigma }}w_{i}\left( R\right) \right) \right] ^{1/2}}%
\int_{R}^{r}t^{1-N}dt \\ 
& +\left( \int_{R}^{r}t^{1+\varepsilon }\overset{d}{\underset{i=1}{\Sigma }}%
p_{i}\left( t\right) dt+\frac{1}{\varepsilon R^{\varepsilon }}\right) \\ 
& \leq C_{R}+\int_{R}^{r}t^{1+\varepsilon }\overset{d}{\underset{i=1}{\Sigma 
}}p_{i}\left( t\right) dt\text{,}%
\end{array}
\label{final}
\end{equation}%
where 
\begin{equation*}
C_{R}=\sqrt[2]{C}\left[ F\left( \overset{d}{\underset{i=1}{\Sigma }}%
w_{i}\left( R\right) \right) \right] ^{-1/2}\frac{R^{N-2}}{N-2}+\frac{1}{%
\varepsilon R^{\varepsilon }}.
\end{equation*}%
By taking $r\rightarrow \infty $ in (\ref{final}) we obtain that $p_{j}$ ($%
j=1,...,d$) satisfies (\ref{13}). This completes the proof of the theorem.

We conclude this paper with some remarks:

\begin{remark}
The time-independent Schr\"{o}dinger equation in quantum mechanics is 
\begin{equation*}
(h^{2}/2m)\Delta u=(V-E)u
\end{equation*}%
where $h=6.625\cdot 10^{-27}$ erg sec is the Planck constant, $m$ is the
mass of a particle moving under the action of a force field described by the
potential $V(x,y,z)$ whose wave function is $u(x,y,z,t)$ and the quantity E is the total energy of the particle, problems which falls into the class of equations discussed here.
\end{remark}

\begin{remark}
If (C1)-(C3) are satisfied then%
\begin{equation*}
\int_{1}^{\infty }\left( \int_{0}^{t}f_{i}\left( s,...,s\right) \right)
^{-1/2}dsdt=\infty \text{ for all }i=1,..,d.
\end{equation*}
\end{remark}

\begin{remark}
(see \cite{CD}) If C1)-C2) and 
\begin{equation*}
\int_{1}^{\infty }\left( \underset{i=1}{\overset{d}{\Sigma }}f_{i}\left(
s,...,s\right) \right) ^{-1}ds=\infty ,
\end{equation*}%
are satisfied, then%
\begin{equation*}
\int_{1}^{\infty }\left( \int_{0}^{t}\underset{i=1}{\overset{d}{\Sigma }}%
f_{i}\left( s,...,s\right) \right) ^{-1/2}dsdt=\infty .
\end{equation*}
\end{remark}

\textbf{Acknowledgement. }The author would like to thank to the editors and
reviewers for valuable comments and suggestions which contributed to improve
this article.

\begin{flushright}
\begin{tabular}{l}
nuna adresa \\ 
Drago\c{s}-P\u{a}tru Covei$^{1,2}$ \\ 
$^{1}$Constantin Brâncu\c{s}i University of Târgu-Jiu, \\ 
Calea Eroilor, No 30, Târgu-Jiu, Gorj, \\ 
România. \\ 
$^{2}$West University of Timi\c{s}oara, \\ 
Bld. Pârvan, No. 4, 300223, Timi\c{s}oara, Timi\c{s}, \\ 
România. \\ 
{e-mail: coveid@yahoo.com}%
\end{tabular}
\end{flushright}

\begin{center}
(Ricevita la 13-an de junio, 2010)
\end{center}

\end{document}